\documentclass[11pt]{article}

\usepackage[T1]{fontenc}
\usepackage{lmodern}
\usepackage{amsmath,amssymb,amsthm}
\usepackage[letterpaper,margin=1in]{geometry}
\usepackage[hidelinks]{hyperref}

\hypersetup{
  pdftitle={The Maximum Diagonal Product of a Stochastic Matrix of Bounded Rank},
  pdfauthor={Yair Lavi}
}

\newtheorem{theorem}{Theorem}[section]
\newtheorem{proposition}[theorem]{Proposition}
\newtheorem{lemma}[theorem]{Lemma}
\newtheorem{corollary}[theorem]{Corollary}
\theoremstyle{remark}
\newtheorem*{remark}{Remark}

\title{The Maximum Diagonal Product of a Stochastic Matrix of Bounded Rank}
\author{Yair Lavi}
\date{21 October 2026}

\begin{document}
\maketitle

\begin{abstract}
Let $A=(a_{ij})$ be a nonnegative row-stochastic or column-stochastic
$n\times n$ matrix with $\operatorname{rank}A\le k$, where $1\le k\le n$.
Write $n=qk+s$, where $0\le s<k$. The author conjectured in 2018 that, for
every permutation $\sigma$ of $\{1,\ldots,n\}$,

\[
\prod_{i=1}^n a_{i,\sigma(i)}
\le q^{-q(k-s)}(q+1)^{-(q+1)s}.
\]

We prove this bound and classify all equality cases.

\end{abstract}



\section{Introduction and main results}

Let $S_n$ be the set of permutations of $[n]=\{1,\ldots,n\}$. For a square
matrix $A=(a_{ij})$ and $\sigma\in S_n$, write
\[
d_\sigma(A)=\prod_{i=1}^n a_{i,\sigma(i)}.
\tag{1.1}
\]

A nonnegative matrix is \textbf{row stochastic}, \textbf{column stochastic}, or
\textbf{doubly stochastic} when all its row sums, all its column sums, or both,
respectively, are one.

Throughout, \(\operatorname{rank}\) denotes ordinary matrix rank.

For a positive integer \(t\), let \(J_t\) be the \(t\times t\) matrix with
every entry equal to \(1/t\).

Given \(1\leq k\leq n\), write
\[
 n=qk+s,\qquad 0\leq s<k,                                  \tag{1.2}
\]
and let \(\rho_{n,k}\) denote the balanced \(k\)-part partition consisting
of \(k-s\) parts \(q\) and \(s\) parts \(q+1\). For a partition
\(\rho=(t_1,\ldots,t_k)\), put

\[
\begin{aligned}
J_\rho&=J_{t_1}\oplus\cdots\oplus J_{t_k},
\qquad J_{n,k}=J_{\rho_{n,k}},\\
D_{n,k}&=\prod_{t\in\rho_{n,k}}t^{-t}
=q^{-q(k-s)}(q+1)^{-(q+1)s}.
\end{aligned}
\]

Here \(\oplus\) denotes direct sum; \(J_{n,k}\) is doubly stochastic and has
rank \(k\). Define
\((\Pi_\sigma)_{ij}=1\) when \(j=\sigma(i)\), and zero otherwise.

When \(A\) is nonnegative and row or column stochastic, \(d_\sigma(A)\leq1\),
with equality only for the full-rank permutation matrix \(\Pi_\sigma\).
In 2018, the author conjectured
\cite[Conjectures 4--6]{lavi2018} that every nonnegative row-stochastic or column-stochastic \(A\)
of rank at most \(k\) satisfies \(d_\sigma(A)\le D_{n,k}\), with equality
exactly when \(A=P^{\mathsf T}J_{n,k}P\Pi_\sigma\) for a permutation matrix
\(P\). We prove this and its homogeneous form.

\begin{theorem}
Let $1\le k\le n$, let $\sigma\in S_n$, and let $A$ be a
nonnegative $n\times n$ matrix with $\operatorname{rank}A\le k$.

\begin{enumerate}
\item If $A$ is row or column stochastic, then

\[
\boxed{d_\sigma(A)\le D_{n,k}.}
\tag{1.3}
\]
Equality holds exactly when
\[
A=P^{\mathsf T}J_{n,k}P\Pi_\sigma
\tag{1.4}
\]
for a permutation matrix $P$. Every matrix in (1.4) has rank $k$.

\item For arbitrary nonnegative $A$, let $r_i$ and $c_j$ be its row and column
sums, and put

\[
R(A)=\prod_{i=1}^n r_i,\qquad C(A)=\prod_{j=1}^n c_j.
\tag{1.5}
\]
Then
\[
\boxed{d_\sigma(A)\le D_{n,k}\min\{R(A),C(A)\}.}
\tag{1.6}
\]

A zero row or column gives equality. If every row and column sum is
positive, define

\[
A^r=\operatorname{diag}(r_1^{-1},\ldots,r_n^{-1})A,
\qquad
A^c=A\operatorname{diag}(c_1^{-1},\ldots,c_n^{-1}).
\tag{1.7}
\]

When $R(A)<C(A)$, equality in (1.6) holds exactly when $A^r$ has form
(1.4). When $C(A)<R(A)$, it holds exactly when $A^c$ has that form. When
$R(A)=C(A)>0$, equality holds exactly when both normalized matrices have
that form; either condition then implies the other.
\end{enumerate}
\end{theorem}

Two specializations have short derivations from classical matrix results.
We state them first; Section 6 proves their bounds and equality cases
without invoking Theorem 1.1.

\textbf{Singular specialization.} When $k=n-1$, the balanced block sizes are
$(2,1,\ldots,1)$. The value bound below follows from Ostrowski's
pairwise refinement of Hadamard's diagonal-dominance criterion \cite{ostrowski1937},
equivalently from Brauer's eigenvalue inclusion theorem \cite{brauer1947}; see Section 6.1.

\begin{corollary}
If $n\ge2$ and $A$ is a singular nonnegative
row-stochastic or column-stochastic $n\times n$ matrix, then for every
$\sigma\in S_n$,

\[
d_\sigma(A)\le\frac14.
\tag{1.8}
\]
Equality holds exactly when
$A=P^{\mathsf T}(J_2\oplus I_{n-2})P\Pi_\sigma$ for a permutation matrix $P$.
When $n=2$, the $I_0$ summand is omitted.
\end{corollary}

\textbf{Divisible specialization.} When $n=qk$, all balanced blocks have size
$q$. The value bound below follows from the trace inequality in
Theorem 1.3 of \cite{fanhoffman1953}, derived there from the Hadamard--Gershgorin
eigenvalue bound, and the arithmetic--geometric mean inequality.
Section 6.2 gives the equality argument without invoking Theorem 1.1.

\begin{corollary}
If $n=qk$ with $q,k\ge1$ and $A$ is a nonnegative
row-stochastic or column-stochastic $n\times n$ matrix with
$\operatorname{rank}A\le k$, then
for every $\sigma\in S_n$,

\[
d_\sigma(A)\le q^{-n}.
\tag{1.9}
\]
Equality holds exactly when
$A=P^{\mathsf T}(J_q^{\oplus k})P\Pi_\sigma$ for a permutation matrix $P$.
\end{corollary}

\textbf{Remaining ranks.} For $k\le n-2$ with $k\nmid n$, the trace estimate
gives only $d_\sigma(A)\le(k/n)^n$, whereas $D_{n,k}<(k/n)^n$ by strict
convexity of $t\log t$. We found no earlier proof of the sharp
value or equality classification in this range. Its upper-half portion
has a particularly simple closed form:

\begin{corollary}
If $n/2<k\le n-2$ and $A$ is a nonnegative
row-stochastic or column-stochastic $n\times n$ matrix with
$\operatorname{rank}A\le k$, then
for every $\sigma\in S_n$,

\[
d_\sigma(A)\le 4^{-(n-k)}.
\tag{1.10}
\]
Equality holds exactly when
$A=P^{\mathsf T}(J_2^{\oplus(n-k)}\oplus I_{2k-n})P\Pi_\sigma$ for a
permutation matrix $P$.
\end{corollary}

\begin{remark}
Each corollary has a homogeneous form: multiply the right-hand
side of its bound by $\min\{R(A),C(A)\}$. Its equality cases, including zero
rows or columns, are those specified in Theorem 1.1(2).
\end{remark}

We note that the bound $d_\sigma(A)\le\operatorname{per}A$, combined with the sharp
permanent bound in \cite{lavi2026}, does not recover (1.3). Indeed, the resulting upper
bound is $D_{n,k}(q!)^{k-s}((q+1)!)^s>D_{n,k}$ whenever $k<n$.

The proof starts by selecting, for each column, a row containing its
maximum, then normalizing the resulting matrix by columns. This does not
lower the diagonal product or raise the rank. Each diagonal entry of the
resulting column-stochastic matrix caps its column. Column filling bounds
the squared Frobenius norm above in terms of these entries, while a
singular-value estimate bounds it below using the trace and rank. A
logarithmic inequality sharp at $1/q$ and $1/(q+1)$ then gives (1.3);
tracing equality through the same steps identifies every equality case.
Section 5 gives the homogeneous form and an exact deficit identity.
Section 6 derives the singular and divisible specializations from classical
results without invoking Theorem 1.1.

\section{Peak rows and the rank constraint}

We begin with a row-stochastic matrix $A$ whose main diagonal is positive.
The assumption loses nothing for an upper bound on its diagonal product: a
zero diagonal entry makes that product zero. Write
$d(A)=\prod_i a_{ii}$.

\begin{lemma}[peak-row reduction]
Let $A$ be nonnegative and row
stochastic, with positive main diagonal. There is a column-stochastic matrix
$C$ such that

\[
\operatorname{rank}C\le\operatorname{rank}A,
\qquad c_{ij}\le c_{jj}\quad(1\le i,j\le n),
\qquad d(C)\ge d(A).
\tag{2.1}
\]
\end{lemma}

\begin{proof}
For each column $j$, choose an original row $i(j)$ attaining
$m_j=\max_i a_{ij}$, and let row $j$ of $B$ be row $i(j)$ of $A$. Repeated
choices are allowed. Then $B$ is row stochastic, its rows lie in the row
space of $A$, and

\[
b_{jj}=m_j\ge a_{jj}>0,
\qquad b_{ij}\le m_j\quad\text{for every }i.
\tag{2.2}
\]

Let $\gamma_j=\sum_i b_{ij}$. These sums are positive, and
$\sum_j\gamma_j=n$ because the rows of $B$ sum to one. The
arithmetic--geometric mean inequality therefore gives
$\prod_j\gamma_j\le1$, with equality only if every $\gamma_j=1$. Set

\[
C=B\operatorname{diag}(\gamma_1^{-1},\ldots,\gamma_n^{-1}).
\tag{2.3}
\]

It is column stochastic. Dividing each column by one positive number
preserves the inequalities in (2.2), while the invertible diagonal factor
preserves rank. Finally,

\[
d(C)=\frac{\prod_jm_j}{\prod_j\gamma_j}
\ge\prod_jm_j\ge d(A).
\tag{2.4}
\]

This proves the lemma.
\end{proof}

The next estimate applies to any column-stochastic matrix with the cap in
(2.1); it does not depend on how that matrix was obtained. For $0<x\le1$,
define

\[
t(x)=\left\lfloor\frac1x\right\rfloor,
\qquad r(x)=1-t(x)x,
\qquad \psi(x)=r(x)\bigl(x-r(x)\bigr),
\qquad h(x)=x+\psi(x).
\tag{2.5}
\]

Here $0\le r(x)<x$. At a reciprocal integer, $r(x)=\psi(x)=0$ and
$h(x)=x$.

\begin{lemma}[column filling]
If $y_1,\ldots,y_n\in[0,x]$ sum to one,
then

\[
\sum_{i=1}^n y_i^2
\le t(x)x^2+r(x)^2=x-\psi(x).
\tag{2.6}
\]
\end{lemma}

\begin{proof}
The feasible vectors form a compact polytope, and the convex
function $\sum_i y_i^2$ attains a maximum at one of its vertices. At a
vertex, at most one coordinate can lie strictly between $0$ and $x$:
otherwise a small amount could be moved in either direction between two
such coordinates. Thus the nonzero coordinates consist of $t(x)$ copies
of $x$ and, when $r(x)>0$, one copy of $r(x)$. Their squared sum is
$t(x)x^2+r(x)^2=x-r(x)(x-r(x))$.
\end{proof}

For a real matrix $M$, write $\|M\|_F^2=\sum_{i,j}m_{ij}^2$. Its nuclear
norm $\|M\|_*$ is the sum of its nonzero singular values.

\begin{proposition}[rank-sensitive diagonal constraint]
Let $C$ be a
nonnegative column-stochastic $n\times n$ matrix of rank at most $k$, with
$c_{ij}\le c_{jj}$ for every $i,j$. Put $x_j=c_{jj}$. Then
\[
\boxed{\sum_{j=1}^n h(x_j)\le k.}
\tag{2.7}
\]
\end{proposition}

\begin{proof}
Each column sums to one and is capped by $x_j$, so
$1/n\le x_j\le1$. Lemma 2.2, summed over the columns, gives

\[
\|C\|_F^2\le\sum_j\bigl(x_j-\psi(x_j)\bigr).
\tag{2.8}
\]

Let $\ell=\operatorname{rank}C\le k$ and let $s_1,\ldots,s_\ell$ be its
nonzero singular values. Since $(s_\alpha-1)^2\ge0$ and
$\operatorname{tr}C\le\|C\|_*$,

\nopagebreak[4]
\[
\|C\|_F^2=\sum_{\alpha=1}^\ell s_\alpha^2
\ge2\sum_{\alpha=1}^\ell s_\alpha-\ell
\ge2\operatorname{tr}C-k.
\tag{2.9}
\]

For the trace inequality, take a singular-value decomposition of $C$:
if $u_\alpha,v_\alpha$ are paired unit singular vectors, then
$\operatorname{tr}C=\sum_\alpha s_\alpha
\langle u_\alpha,v_\alpha\rangle\le\sum_\alpha s_\alpha$.
Since $\operatorname{tr}C=\sum_jx_j$, combining (2.8) and (2.9) yields
(2.7).
\end{proof}

\begin{remark}
Fan and Hoffman's rank inequality (Theorem 1.2 of \cite{fanhoffman1953}), applied to
$C^{\mathsf T}$, and Lemma 2.2 give the stronger estimate
$\sum_j x_j^2/(x_j-\psi(x_j))\le k$. Indeed, the difference between each
summand and $h(x_j)$ is
$\psi(x_j)^2/(x_j-\psi(x_j))\ge0$. We use (2.9) because its
singular-value terms track equality in Section 4.
\end{remark}

The column normalization in Lemma 2.1 is essential to this argument: the
column-filling bound uses columns whose sums are exactly one.

\section{A logarithmic inequality with two contacts}

Fix the quotient $q$ in (1.2), and set

\[
\lambda_q=q(q+1)\log\frac{q+1}{q},
\qquad
\beta_q=-\log q-\frac{\lambda_q}{q}.
\tag{3.1}
\]

\begin{lemma}
For every $0<x\le1$,

\[
\boxed{\log x\le\lambda_q h(x)+\beta_q,}
\tag{3.2}
\]
with equality exactly at $x=1/q$ and $x=1/(q+1)$.
\end{lemma}

\begin{proof}
For $t\ge1$ and $1/(t+1)<x\le1/t$, (2.5) gives

\[
h(x)=2(t+1)x-1-t(t+1)x^2.
\tag{3.3}
\]

The formula also holds at the left endpoint, where both sides equal $x$. Let
$G_q(x)=\lambda_qh(x)+\beta_q-\log x$. At a reciprocal endpoint,

\[
G_q(1/t)=F_q(t)+\beta_q,
\qquad F_q(u)=\log u+\frac{\lambda_q}{u}.
\tag{3.4}
\]

The strict logarithm bounds
$z/(1+z)<\log(1+z)<z$, applied at $z=1/q$, give
$q<\lambda_q<q+1$. Moreover,

\[
F_q'(u)=\frac{u-\lambda_q}{u^2},
\qquad F_q(q)=F_q(q+1).
\tag{3.5}
\]

Hence the smallest value of $F_q(t)$ over positive integers occurs
exactly at $t=q,q+1$. The choice of $\beta_q$ makes $G_q$ zero at both
corresponding reciprocal endpoints and positive at every other reciprocal
endpoint.

It remains to rule out a smaller value inside an interval. Differentiating
(3.3) gives

\[
xG_q'(x)=2\lambda_q(t+1)x(1-tx)-1,
\qquad
\frac{d}{dx}\bigl[xG_q'(x)\bigr]
=2\lambda_q(t+1)(1-2tx).
\tag{3.6}
\]

For $1/(t+1)<x<1/t$, the last quantity is strictly negative: indeed
$x>1/(t+1)\ge1/(2t)$. Thus $xG_q'(x)$ decreases strictly, and $G_q'$ can
change sign only from positive to negative. There is no interior minimum;
every interior value is strictly above the smaller endpoint value. The
endpoint calculation proves (3.2) and its equality statement.
\end{proof}

\section{Proof of the stochastic theorem}

We first prove Theorem 1.1(1) for the main diagonal of a row-stochastic
matrix. If $d(A)=0$, the upper bound is immediate and equality is
impossible because $D_{n,k}>0$. Otherwise apply Lemma 2.1 to obtain $C$,
and write $x_j=c_{jj}$. Lemma 3.1 and Proposition 2.3 give

\[
\begin{aligned}
\log d(A)
&\le\log d(C)=\sum_j\log x_j\\
&\le\lambda_q\sum_j h(x_j)+n\beta_q\\
&\le\lambda_q k+n\beta_q
=\log D_{n,k}.
\end{aligned}
\tag{4.1}
\]

For the last identity, use $n=qk+s$ and
$\lambda_q/q=(q+1)\log((q+1)/q)$. This proves the value.

We now determine equality. Use the construction in the proof of Lemma 2.1,
with any choice of maximizing rows, and retain the resulting
$m_j,B,\gamma_j,C$. If $d(A)=D_{n,k}$, every inequality in
(2.4) and (4.1) is an equality. In particular,
$a_{jj}=m_j=b_{jj}$ for every $j$, and
$\prod_j\gamma_j=1$. Equality in arithmetic--geometric mean, with
$\sum_j\gamma_j=n$, gives $\gamma_j=1$ for every $j$. Thus $B=C$ is doubly
stochastic.

Equality in the scalar inequality (3.2) says that every diagonal entry
$x_j$ of $B$ lies in $\{1/q,1/(q+1)\}$. Hence $\psi(x_j)=0$ and
$h(x_j)=x_j$. Since $\lambda_q>0$, equality in the last inequality of
(4.1) gives $\sum_jh(x_j)=k$. Thus
$\operatorname{tr}B=\sum_jx_j=k$. The two bounds
(2.8)--(2.9) now squeeze $\|B\|_F^2$ to $k$. More precisely,

\[
\|B\|_F^2-2\operatorname{tr}B+k
=\sum_{\alpha=1}^\ell(s_\alpha-1)^2
  +2\bigl(\|B\|_*-\operatorname{tr}B\bigr)
  +(k-\ell)=0,
\tag{4.2}
\]
where $\ell=\operatorname{rank}B$. All three terms on the right are
nonnegative. Consequently $\ell=k$, every nonzero singular value of $B$
equals one, and $\|B\|_*=\operatorname{tr}B=k$.

Since these singular values are all one, write a thin singular-value
decomposition as $B=UV^{\mathsf T}$, where the
$n\times k$ matrices $U,V$ have orthonormal columns. Then
\[
\|U-V\|_F^2
=2k-2\operatorname{tr}(V^{\mathsf T}U)
=2k-2\operatorname{tr}B=0.
\]

Thus $U=V$ and $B=UU^{\mathsf T}$ is a symmetric idempotent. Its entries
are nonnegative and its diagonal is positive. The relation $b_{ij}>0$ is
an equivalence relation on $[n]$: symmetry and reflexivity are clear, and

\[
b_{il}=(B^2)_{il}
=\sum_v b_{iv}b_{vl}\ge b_{ij}b_{jl}>0
\]
proves transitivity. Accordingly, $B$ is a direct sum of strictly positive
stochastic blocks after a simultaneous permutation of rows and columns.
Each such block is idempotent. Its fixed vectors are constant: at a
coordinate where a fixed vector is maximal, strict positivity of all
weights in that row forces every coordinate to have the same value.
Because an idempotent's image is its fixed space, the block has rank one.
Being symmetric and stochastic, it is $J_t$ for its block size $t$.

Every diagonal entry in a $J_t$ block is $1/t$, so the contacts from
Lemma 3.1 force $t\in\{q,q+1\}$. There are $\operatorname{rank}B=k$
blocks. Their sizes sum to $n=qk+s$, hence exactly $k-s$ blocks have size
$q$ and $s$ have size $q+1$. We have shown that
$B=P^{\mathsf T}J_{n,k}P$.

It remains to recover $A$; equality of its diagonal with $B$ does not by
itself say that the matrices are equal. The rows of $B$ are rows of $A$,
and both matrices have rank $k$, since
$k=\operatorname{rank}B\le\operatorname{rank}A\le k$. Their row spaces
therefore coincide. Every row of $A$ is constant on each block of $B$.
If row $j$ belongs to a block of size $t$, its diagonal entry is
$a_{jj}=1/t$, so all $t$ entries of that row in the same block equal
$1/t$. They already sum to one; nonnegativity and row stochasticity force
all other entries to vanish. Row $j$ of $A$ is row $j$ of $B$, and hence
$A=B$.

Conversely, $P^{\mathsf T}J_{n,k}P$ is row stochastic, has rank $k$, and
its diagonal product is $D_{n,k}$. This completes the main-diagonal,
row-stochastic case. For a prescribed $\sigma$, first replace $A$ by
$A\Pi_\sigma^{\mathsf T}$, and then transpose if the result is column
stochastic. The resulting matrix is row stochastic, has the same rank,
and has main-diagonal product $d_\sigma(A)$. At equality it is
$P^{\mathsf T}J_{n,k}P$, which is symmetric. Hence
$A\Pi_\sigma^{\mathsf T}=P^{\mathsf T}J_{n,k}P$ even if a transpose was used,
so $A=P^{\mathsf T}J_{n,k}P\Pi_\sigma$. Conversely, every matrix of this
form has rank $k$ and attains $D_{n,k}$. This proves Theorem 1.1(1).
\qed

\section{Homogeneous form and an exact deficit}

\begin{proof}[Proof of Theorem 1.1(2)]
A zero row or column makes both sides of (1.6)
zero. Suppose every row and column sum is positive. Positive row and column scaling
preserve rank, while

\[
d_\sigma(A)=R(A)d_\sigma(A^r)
=C(A)d_\sigma(A^c).
\tag{5.1}
\]

Theorem 1.1(1) applies to the row-stochastic matrix $A^r$ and the
column-stochastic matrix $A^c$. Taking the smaller of the resulting bounds
proves (1.6). If one of the two products is strictly smaller, equality is
equivalent to equality for its corresponding normalization. If the two
products agree, equality in either bound implies equality in both, since
the two expressions in (5.1) have the same numerator. The equality cases
are therefore exactly those stated in Theorem 1.1(2).
\end{proof}

The proof also separates the ways in which a positive diagonal product
can fall below the sharp value. The following identity is stated for a
row-stochastic matrix and the main diagonal; Theorem 1.1 transports the
bound to the other cases.

\begin{proposition}[deficit identity]
Let $A$ be row stochastic with
$d(A)>0$ and $\operatorname{rank}A\le k$. For any choice of maximizing rows
in the proof of Lemma 2.1, use the resulting $m_j,B,\gamma_j,C$ and put
$x_j=c_{jj}$. Let $\ell=\operatorname{rank}C$ and let
$s_1,\ldots,s_\ell$ be the nonzero singular values of $C$. Define

\[
G_q(x)=\lambda_qh(x)+\beta_q-\log x.
\]
Then
\[
\begin{aligned}
\log\frac{D_{n,k}}{d(A)}
={}&\sum_j\log\frac{m_j}{a_{jj}}
    -\sum_j\log\gamma_j+\sum_jG_q(x_j)\\
 &+\lambda_q\left[
     \sum_{\alpha=1}^\ell(s_\alpha-1)^2
     +2\bigl(\|C\|_*-\operatorname{tr}C\bigr)
     +(k-\ell)
     +\sum_j\left(x_j-\psi(x_j)-\sum_i c_{ij}^2\right)
   \right].
\end{aligned}
\tag{5.2}
\]

Each of the four groups on the right of (5.2) is nonnegative, as is each
term inside the brackets.
\end{proposition}

\begin{proof}
Lemma 2.1 gives $m_j\ge a_{jj}$ and
$\prod_j\gamma_j\le1$. Lemma 3.1 gives $G_q(x_j)\ge0$. The terms in
brackets are nonnegative by singular-value squares, the nuclear-norm trace
inequality, $\ell\le k$, and Lemma 2.2, respectively. To verify the
identity, first use (3.1) and (4.1) to write
\[
\log D_{n,k}-\log d(C)
=\sum_jG_q(x_j)
 +\lambda_q\left(k-\sum_jh(x_j)\right).
\]

Since $h(x)=x+\psi(x)$, the expression in parentheses is

\[
\bigl(\|C\|_F^2-2\operatorname{tr}C+k\bigr)
 +\sum_j\left(x_j-\psi(x_j)-\sum_i c_{ij}^2\right).
\]

Expand the first term by the identity in (4.2), with $C$ in place of $B$,
and use $d(C)=\prod_jm_j/\prod_j\gamma_j$. This gives (5.2).
\end{proof}

\section{Classical proofs in two specializations}

We give those proofs, including the equality cases.
In each case, replace $A$ by $A\Pi_\sigma^{\mathsf T}$ and transpose if
necessary. The result is row stochastic, has the same rank and diagonal
product, and has the chosen product on its main diagonal. The equality
matrices obtained below are symmetric, so reversing these operations gives
the forms stated in the corollaries.

\subsection{Singular specialization}

Let $B$ be a singular nonnegative row-stochastic matrix. Ostrowski's
nonsingularity criterion \cite{ostrowski1937}, or the equivalent eigenvalue theorem of
Brauer \cite{brauer1947} applied at zero, gives distinct $i,j$ such that

\[
b_{ii}b_{jj}
\le \Bigl(\sum_{u\ne i}b_{iu}\Bigr)
     \Bigl(\sum_{u\ne j}b_{ju}\Bigr)
=(1-b_{ii})(1-b_{jj}).
\]

Thus $b_{ii}+b_{jj}\le1$. Every diagonal entry of $B$ lies in $[0,1]$, so

\[
d(B)\le b_{ii}b_{jj}
\le\left(\frac{b_{ii}+b_{jj}}2\right)^2\le\frac14.
\]

At equality, $b_{ii}=b_{jj}=1/2$ and every other diagonal entry is one.
Every row outside $\{i,j\}$ is therefore the corresponding unit row.
After a simultaneous permutation putting $i,j$ first, singularity gives
$0=\det B=1/4-b_{ij}b_{ji}$. Row stochasticity gives
$b_{ij},b_{ji}\le1/2$, so both equal $1/2$; those two rows have no other
nonzero entries. Hence $B=J_2\oplus I_{n-2}$ after the permutation.
Conversely, this matrix is singular and has diagonal product $1/4$.

\subsection{Divisible specialization}

Suppose $n=qk$ and let $B$ be a nonnegative row-stochastic matrix of rank
at most $k$. Theorem 1.3 of \cite{fanhoffman1953}, proved there
using the Hadamard--Gershgorin eigenvalue bound, gives
$\operatorname{tr}B\le\operatorname{rank}B\le k$. The
arithmetic--geometric mean inequality therefore gives

\[
d(B)\le\left(\frac{\operatorname{tr}B}{n}\right)^n
\le\left(\frac{k}{n}\right)^n=q^{-n}.
\]

If equality holds, every diagonal entry is $1/q$ and
$\operatorname{tr}B=\operatorname{rank}B=k$. Let $m$ count the nonzero
eigenvalues of $B$ with algebraic multiplicity. All eigenvalues lie in
the closed unit disk, while $m\le\operatorname{rank}B=k$. Taking real
parts in their sum gives
$k=\operatorname{tr}B\le m\le k$; hence $m=k$ and every nonzero
eigenvalue is one. The zero eigenvalue is semisimple because its
algebraic multiplicity $n-k$ equals $\dim\ker B$. The eigenvalue one
is semisimple because $B^r$ is row stochastic for every positive integer
$r$, so its powers are bounded. Therefore $B^2=B$.

Every diagonal entry of $B$ is positive. Flor's classification of
nonnegative idempotents \cite[Theorem 2]{flor1969}, in notation adapted here, states
that a rank-$k$ nonnegative idempotent is permutation-similar to
\[
\begin{pmatrix}
H & HX & 0 & 0\\
0 & 0 & 0 & 0\\
YH & YHX & 0 & 0\\
0 & 0 & 0 & 0
\end{pmatrix},
\qquad H=H_1\oplus\cdots\oplus H_k,
\]
where $X,Y$ are nonnegative matrices of compatible sizes and each $H_i$
is a nonnegative rank-one idempotent; empty block groups are allowed.
Only $H$ has nonzero diagonal entries, so positivity of every diagonal
entry of $B$ makes the other block groups empty. Write each
$H_i=u_i v_i^{\mathsf T}$ with $u_i,v_i\ge0$. Its positive diagonal forces
all coordinates of $u_i,v_i$ to be positive. Thus $B$ is, after a
simultaneous permutation, a direct sum of positive rank-one stochastic
blocks. Each block has identical rows given by a probability vector $p$;
its diagonal entries are the coordinates of $p$. Since each
diagonal entry is $1/q$, every block has size $q$ and is $J_q$.
There are $k$ blocks because $\operatorname{rank}B=k$. Conversely,
$J_q^{\oplus k}$ has rank $k$ and diagonal product $q^{-n}$.

Finally, (5.1) transfers both specializations and their equality
classifications to arbitrary nonnegative matrices with the same rank
restrictions. A zero row or column gives equality in the homogeneous
bound.

\section*{Acknowledgments}

The proof of this conjecture was carried out by GPT-5.6-sol, GPT-6 Astra,
Claude Fable 5, and Claude Fable 5.1, under the guidance of the author. The author has reviewed the
resulting proof arguments. Responsibility for the final text rests with the
author.

\end{document}